\documentclass[reqno,12pt]{amsart}

\newtheorem{theorem}{Theorem}[section]
\newtheorem{lemma}[theorem]{Lemma}
\newtheorem{corollary}[theorem]{Corollary}

\theoremstyle{definition}

\theoremstyle{remark}
\newtheorem{remark}[theorem]{Remark}

\newcommand{\mysection}[1]{\section{#1}
\setcounter{equation}{0}}

\newcommand{\bC}{\mathbb C}
\newcommand{\bR}{\mathbb R}

\renewcommand{\epsilon}{\varepsilon}

\begin{document}
\title[Schr\"odinger equations
]{Unique continuation for the Schr\"odinger equation with  gradient
vector potentials}

\author[H. Dong]{Hongjie Dong}
\address[H. Dong]
{Department of Mathematics, University of Chicago, 5734 S.
University Avenue,  Chicago, Illinois 60637, USA}
\email{hjdong@math.uchicago.edu}

\author[W. Staubach]{Wolfgang Staubach}
\address[W. Staubach]
{Department of Mathematics, University of Chicago, 5734 S.
University Avenue,  Chicago, Illinois 60637, USA}
\email{wolf@math.uchicago.edu}



\begin{abstract}
We obtain unique continuation results for Schr\"odinger equations
with time dependent gradient vector potentials. This result with an
appropriate modification also yields unique continuation properties
for solutions of certain nonlinear Schr\"odinger equations.
\end{abstract}

\maketitle

\mysection{Introduction}

In this paper we study unique continuation properties of solutions
to the Schr\"odinger equations having time dependant gradient
magnetic vector potentials. Specifically, we consider equations of
the form
$$
i\partial_t u+\Delta_{x} u= p_1(x,t)\cdot \nabla_{x}u+p_2(x,t)\cdot
\nabla_{x}\bar u
$$
\begin{equation}
                                   \label{eq10.17}
 + V_1(x,t)u+V_2(x,t)\bar
u \quad \text{on}\,\,\, \mathbb{R}^{d} \times (0,1).
\end{equation}
where the vector potentials $p_1$ and $p_2$ and the scalar
potentials $V_1,V_2$ are assumed to belong to certain Banach spaces
considered in the work of Ionescu and Kenig \cite{IonKen2}. The
unique continuation results for solutions to equation
\eqref{eq10.17} can be used to obtain unique continuation properties
of solutions to nonlinear Schr\"odinger equations of type
\begin{equation}
                                   \label{eq10.18}
i\partial_t u+\Delta_{x} u= P(u,\bar u,\nabla_x u,\nabla_x \bar u)
\,\,\, \text{on}\,\,\, \mathbb{R}^{d} \times (0,1),
\end{equation}
where $P:\bC\times \bC\times \bC^d\times \bC^d\to \bC$ is a
polynomial.

Unique continuation for the solutions of
 Schr\"odinger equations has been studied by
various authors e.g. Bourgain \cite {B}, Kenig and Ionescu
\cite{Ionken1}, \cite{IonKen2}, Kenig, Ponce and Vega \cite{KPV},
and Escauriaza, Kenig, Ponce and Vega \cite{EKPV2}.

The main goal here is to find sufficient conditions on the solution
$u$ of the linear equation \eqref{eq10.17}, at two different times
$t=0$ and $t=1$, which guarantee that $u\equiv 0$. In the case of
nonlinear Schr\"odinger equation \eqref{eq10.18} the uniqueness
result is obtained from the information on the difference of two
possible solutions at two different times.

Our approach follows that of \cite{EKPV2} but since we also have to
deal with the presence of the gradient vector potential, we need
some estimates in the spirit of \cite{IonKen2}.

Let $H=i\partial_t+\Delta$ denote the Schr\"odinger operator and
denote by $H^1$ the Sobolev space of tempered distributions $u$ for
which $(1-\Delta)^{\frac{1}{2}}u\in L^{2}(\mathbb{R}^{d})$. For any
$k\in \mathbb{Z}^{d}$ let $Q_{k}$ denote the cube
$$\{x\in \mathbb{R}^{d}; x_j\in [k_j -\frac{1}{2},
k_j+\frac{1}{2}), \,\,\, j=1,\dots,\,d\}.$$ Following
\cite{IonKen2}, for $p$, $q\in [1,\infty]$ we define the Banach
space $B^{p,q}(\mathbb{R}^d)$ using the norm
\begin{equation*}
\Vert f\Vert_{B^{p,q}}:=\Big\{\sum_{k\in \mathbb{Z}^{d}}\Vert
f\Vert^{p}_{L^{q}(Q_{k})}\Big\}^{\frac{1}{p}},\,\,\, 1\leq p<\infty
,\,\,\, \Vert f\Vert_{B^{\infty,q}}:= \sup_{k\in
\mathbb{Z}^{d}}\Vert f\Vert_{L^{q}(Q_{k})}
\end{equation*}
From the definition it follows that $B^{p,p}=L^{p}$, $1\leq p\leq
\infty$, and
\begin{equation*}
B^{p_1 , q_1}\hookrightarrow B^{p_2 , q_2}\,\,\,\text{if}\,\,\,
q_1\geq q_2\,\,\, \text{and}\,\,\, p_1 \leq p_2 .
\end{equation*}
Our main result for the linear Schr\"odinger equation is the
 following:
\begin{theorem}
                                    \label{thm1}
Let $u\in C([0,1]:H^1)$ be a solution of the equation
\begin{equation}
                                    \label{eq10.21}
Hu=V_1u+V_2\bar u+p_1\cdot\nabla_x u+p_2\cdot\nabla_x \bar
u\quad\text{on}\,\,\bR^d\times(0,1),
\end{equation}
where
\begin{equation}
                                      \label{eq10.24}
V_1,V_2\in B_x^{2,\infty}L_t^{\infty}(\bR^d\times[0,1]),\quad
|p_1|,|p_2|\in B_x^{1,\infty}L_t^{\infty}(\bR^d\times[0,1]).
\end{equation}
If there exists $a\geq c_0(d,\|p_1\|_{L_{t,x}^{\infty}},
\|p_2\|_{L_{t,x}^{\infty}})$ such that
\begin{equation}
                                        \label{eq10.49}
u(\cdot,0),\,u(\cdot,1)\in H^1(e^{a|x|^2}),
\end{equation}
then $u\equiv 0$.
\end{theorem}

\begin{remark}
In some sense, this theorem is a generalization of the result
obtained in \cite{EKPV2}, to the case of Schr\"odinger operators
with gradient vector potentials.  To satisfy \eqref{eq10.24}, we
need $p_j,V_j,j=1,2$ to be bounded  and decay with certain rates at infinity. For example, all Schwarz
 functions and all functions of the form
$1/(1+|x|^{d+\alpha}),\alpha>0$ satisfy this assumption. By going
through the proof of Theorem \ref{thm1} and using Sobolev embedding
theorem, one can see that the condition of $V$ in \eqref{eq10.24}
can be relaxed to
$$
V_1,V_2\in
\big(B_x^{2,\infty}+B_x^{\frac{2d}{d+1},2d}\big)L_t^{\infty}(\bR^d\times[0,1])\quad
\text{if}\,\,d\geq 2,
$$
$$
V_1,V_2\in\big(B_x^{2,\infty}+
B_x^{\frac{1}{1-\gamma},\frac{2}{1-2\gamma}}\big)L_t^{\infty}(\bR^d\times[0,1])\quad
\text{if}\,\,d=1,
$$
for any $\gamma\in (0,1/2]$. However, at present we don't know if
\eqref{eq10.24} can be improved significantly.
\end{remark}

We show that for the nonlinear Schr\"odinger equation, one has the
following unique continuation property.
\begin{theorem}
                                          \label{thm2}
Consider
\begin{equation}
                                        \label{eq10.37}
Hu=P(u,\bar u,\nabla_x u,\nabla_x \bar u)\quad
\text{on}\,\,\bR^d\times (0,1),
\end{equation}
 and let $X$ be either $L^2(\bR^d)$ or $B^{1,2}(\bR^d)$ and  $J^{s}:=(1-\Delta_{x})^{\frac{s}{2}}$.
 Assume that
$u_1,u_2$ are solutions to \eqref{eq10.37} and we have
$$
J^s u_1, J^s u_2\in C([0,1]:X),\,\,s>d/2+3,
$$
$$
\|J^s u_1\|_{C([0,1]:X)}\leq N,\quad \|J^s u_2\|_{C([0,1]:X)}\leq N,
$$
$$
|P(z_1,z_2,w_1,w_2)|\le N(|z_1|^2+|z_2|^2+|w_1|^2+|w_2|^2),
$$
for a constant $N>0$ and any $|(z_1,z_2,w_1,w_2)|\leq 1$. If
$u=u_1-u_2$ satisfies \eqref{eq10.49} for some $a\geq c_0(d,N)$,
then $u_1\equiv u_2$ in $\bR^d\times [0,1]$.
\end{theorem}

\begin{remark}
                    \label{rem1}
The result of Theorem \ref{thm2} can be extended to more general
nonlinearities of the form $F(x,t,u,\bar u, \nabla_x u,\nabla_x \bar
u)$ where one only needs to assume that
$$
F(x,t,u_1,\bar u_1, \nabla_x u_1,\nabla_x \bar u_1)-F(x,t,u_1,\bar
u_1, \nabla_x u_1,\nabla_x \bar u_1)
$$
\begin{equation*}
=p_1(x,t)\cdot \nabla_{x}u +p_2(x,t)\cdot
\bar\nabla_{x}u+V_1(x,t)u+V_2(x,t)\bar u
\end{equation*}
with $u= u_1-u_2$, for some $V_1,V_2,p_1,p_2$ as in \eqref{eq10.24}
(See Kenig et al. \cite{IonKen2} and \cite{KPV3}).
\end{remark}

The rest of the article is organized as follows: In section 2 we
shall prove suitable Carleman estimates which are the main tools in
obtaining unique continuation results. In section 3 we shall
establish lower bounds for the $L^{2}$ space-time norm of the
solutions to \eqref{eq10.17} and \eqref{eq10.18}. These lower bounds
combined with the Carleman estimates of section 2 yield the desired
unique continuation properties.

\section*{Acknowledgment}
The authors would like to express their gratitude to Carlos Kenig
for bringing this problem to their attention and many enlightening
discussions.

\mysection{Carleman estimates and upper bounds}

Let $\phi$ be a smooth function on $\bR$ with the following
properties: $\phi(0)=0$, $\phi'$ nonincreasing, $\phi'(r)=1$ if
$r\leq 1$ and $\phi'(r)=0$ if $r\geq 2$. For any $\lambda\geq 1$,
let $\phi_\lambda(r)=\lambda\phi(r/\lambda)$. As in \cite{IonKen2},
for any $\beta>0,\lambda>0$, we define the Banach space
$X=X_{\lambda}(\bR^d\times \bR)$ using the norm
$$
\|f\|_X:=\inf_{f_1+f_2=f}[\|f_1\|_{B^{1,2}_xL^2_t}
+\lambda^3\|J^1f_2\|_{L^1_tL_x^2}].
$$
We also define the Banach Space $X'=X'_{\beta
,\lambda}(\mathbb{R}^{d}\times \mathbb{R})$ using the norm
$$
\|u\|_{X'}:=\max\{\|J^{1/2}u\|_{L^\infty_tL^2_x},\|J^1u\|_{B^{\infty,2}_xL^2_t},
\|(1+\beta\phi'_\lambda(x_1))u\|_{B^{\infty,2}_xL^2_t}\}.
$$
We will need the following Carleman type estimate:
\begin{lemma}
                                    \label{lem1}
There are constant $\bar C,C_0$ and $m$ only depending on $d$ such
that
\begin{equation*}
\|e^{\beta\phi_\lambda(x_1)}u\|_{L_x^2L_t^2}+
\|e^{\beta\phi_{\lambda}(x_1)}|\nabla u|\|_{B^{\infty,2}_xL^2_t}
\leq \bar C\|e^{\beta\phi_\lambda(x_1)}Hu\|_X
\end{equation*}
for any $u\in C(\bR:H^1)$ with $u(\cdot,t)\equiv 0$ for $t\notin
[0,1]$, any $\beta\in [1,\infty)$ and any $\lambda\geq C_0\beta^m$.
\end{lemma}
\begin{proof}

In \cite{IonKen2}, Ionescu and Kenig established the Carleman
estimate
$$
\Vert e^{\beta\phi_\lambda(x_1)} u\Vert_{X'} \leq \bar
C\|e^{\beta\phi_\lambda(x_1)}Hu\|_X,
$$
under the same assumptions as in Lemma \ref{lem1}. Using the
$B^{\infty ,2}$ boundedness of pseudodifferential operators of order
zero obtained in \cite{smoothing}, one can show that
\begin{equation*}
 \|e^{\beta\phi_\lambda(x_1)}u\|_{L_x^2L_t^2}+
\|e^{\beta\phi_{\lambda}(x_1)}|\nabla u|\|_{B^{\infty,2}_xL^2_t}
\leq C \max\big\{\|J^{1/2}
(e^{\beta\phi_\lambda(x_1)}u)\|_{L^\infty_tL^2_x},
\end{equation*}
$$
\|J^1(e^{\beta\phi_\lambda(x_1)}u)\|_{B^{\infty,2}_xL^2_t},
\|(1+\beta\phi'_\lambda(x_1))e^{\beta\phi_\lambda(x_1)}u\|_{B^{\infty,2}_xL^2_t}\big\}
$$ for some constant $C$.
\end{proof}

Without loss of generality, in the sequel we assume that $m\geq 2$.
Lemma \ref{lem1} implies the following theorem.

\begin{theorem}

 \label{thm2.2} There exist constant $C,C_0,m\geq 2$ only
depending on $d$ such that
$$
\|e^{\beta\phi_{\lambda}(|x|/\sqrt{d})}v\|_{L^2_xL^2_t}+
\|e^{\beta\phi_{\lambda}(|x|/\sqrt{d})}|\nabla
v|\|_{B^{\infty,2}_xL^2_t}
$$
\begin{equation}
                                \label{eq5.05}
\leq C\|e^{\beta\phi_{\lambda}(|x|)}Hv\|_X+
C\lambda^3\sum_{j=0}^1\|J^1(e^{\beta\phi_{\lambda}(|x|)}v(x,j))\|_{L_x^2}
\end{equation}
for any $v\in C(\bR:H^1)$, $\beta\in [1,\infty)$ and any
$\lambda\geq C_0\beta^m$.
\end{theorem}
\begin{proof}
First we claim that there exist constant $C,m\geq 2$ only depending
on $d$ such that
$$\|e^{\beta\phi_{\lambda}(x_1)}v\|_{L^2_xL^2_t}+
\|e^{\beta\phi_{\lambda}(x_1)}|\nabla v|\|_{B^{\infty,2}_xL^2_t}$$
\begin{equation}
                                \label{eq5.24}
\leq C\|e^{\beta\phi_{\lambda}(x_1)}Hv\|_X+
C\lambda^3\sum_{j=0}^1\|J^1(e^{\beta\phi_{\lambda}(x_1)}v(x,j))\|_{L_x^2}
\end{equation}
for any $v\in C(\bR:H^1)$, $\beta\in [1,\infty)$ and any
$\lambda\geq C_0\beta^m$.

To prove \eqref{eq5.24}, for any $\epsilon\in (0,1/10)$ we define a
smooth cutoff function $\psi_\epsilon:\bR\to [0,1]$ such that
$$\left\{
  \begin{array}{ll}
    \psi_\epsilon\equiv 1, & \text{on}\,\,[2\epsilon,1-2\epsilon];\\
    \psi_\epsilon\equiv 0, & \text{on}\,\,(-\infty,\epsilon]\cup
[1-\epsilon,\infty),
  \end{array}
\right.$$ and $\psi_\epsilon$ is increasing on $(-\infty, 1/2)$ and
decreasing on $(1/2,\infty)$. Since $u(x,t)=v(x,t)\psi_\epsilon(t)$
satisfies the condition in Lemma \ref{lem1}, we have
$$\|e^{\beta\phi_{\lambda}(x_1)}v\psi_\epsilon\|_{L^2_xL^2_t}+
\|e^{\beta\phi_{\lambda}(x_1)}|\psi_\epsilon\nabla
v|\|_{B^{\infty,2}_xL^2_t}$$
$$
\leq C\|e^{\beta\phi_{\lambda}(x_1)}Hv\|_X+
C\lambda^3\|J^1(e^{\beta\phi_{\lambda}(x_1)}\psi_\epsilon'v(x,t))\|_{L_t^1L_x^2}
$$
Letting $\epsilon\to 0$ yields \eqref{eq5.24} immediately. Notice
that this estimate still holds with $x_j$ or $-x_j$ in place of
$x_1$ where $j=1,2,\cdots,d$. Thus to prove \eqref{eq5.05} it
suffices to use the monotonicity of $\phi_{\lambda}$, the obvious
inequalities
\begin{equation}
                                   \label{exp}
 |x|/\sqrt{d}\leq \max_j|x_j|,\quad
e^{\beta\phi_{\lambda}(|x|/\sqrt{d})}\leq
\sum_{j=1}^d(e^{\beta\phi_{\lambda}(x_j)}+e^{\beta\phi_{\lambda}(-x_j)})
\end{equation}
and the triangle inequality. This finishes the proof of the theorem.
\end{proof}

We define the operator
$$
H_{V,p}u=Hu-V_1u-V_2\bar u-p_1\cdot\nabla u-p_2\cdot\nabla \bar u.
$$

\begin{corollary}
                                                \label{cor1}
Under the assumption of Theorem $\ref{thm2.2}$ and the additional
assumptions
\begin{equation}
                                                       \label{extraassumption}
\||p_1|\|_{B^{1,\infty}_xL^\infty_t}+\||p_2|\|_{B^{1,\infty}_xL^\infty_t}\leq
\frac{1}{2C},\quad
\|V_1\|_{B^{2,\infty}_xL^\infty_t}+\|V_2\|_{B^{2,\infty}_xL^\infty_t}\leq
\frac{1}{2C},
\end{equation}
with $C$ as in Theorem $\ref{thm1}$, we have
$$\|e^{\beta\phi_{\lambda}(|x|/\sqrt{d})}v\|_{L^2_xL^2_t}+
\|e^{\beta\phi_{\lambda}(|x|/\sqrt{d})}|\nabla
v|\|_{B^{\infty,2}_xL^2_t}$$
\begin{equation*}
\leq N\|e^{\beta\phi_{\lambda}(|x|)}H_{V,p}v\|_X+
N\lambda^3\sum_{j=0}^1\|J^1(e^{\beta\phi_{\lambda}(|x|)}v(x,j))\|_{L_x^2},
\end{equation*}
with $N=2C$.
\end{corollary}
\begin{proof}
This follows by applying the inequality \eqref{eq5.24} to $Hv=
H_{V,p}v+V_1 v +V_2\bar v +p_1\cdot\nabla v+p_2\cdot\nabla \bar v$
and observing that the definition of the space $X$, the H\"older
inequality and the assumptions \eqref{extraassumption} imply
$$
C\Vert e^{\beta\phi_{\lambda}(x_1)}(V_1 v +V_2\bar v +p_1\cdot\nabla
v+p_2\cdot\nabla \bar v) \Vert_{X}
$$
$$\leq
\frac{1}{2}(\|e^{\beta\phi_{\lambda}(x_1)}v\|_{L^2_xL^2_t}+
\|e^{\beta\phi_{\lambda}(x_1)}|\nabla v|\|_{B^{\infty,2}_xL^2_t}).$$
The result of the corollary follows now by using \eqref{exp}.
\end{proof}
We are now ready to state and prove the main theorem of this
section.
\begin{theorem}
                        \label{thm2.4}
Let $u\in C([0,1]:H^1)$ be a solution of \eqref{eq10.21} and assume
that for some $a>0$ one has
\begin{equation*}
u(\cdot,0),\,u(\cdot,1)\in H^1(e^{a|x|^2}),\quad
\sum_{j=0}^1\|u(x,j)\|_{H^1(e^{a|x|^2})}\leq K,
\end{equation*}
for some constant $K>0$. Then there exist $N=N(d)>0$ and $R_0>0$
depending only on $d,a\wedge 1,K$, $p_1,p_2,V_1,V_2$, the
$L^2([0,1]:H^1)$ norm of $u$, such that for all $R\geq R_{0}$
$$
\|u\|_{L^2(\{R\leq |x|\leq R+1\}\times [0,1])}+\||\nabla
u|\|_{L^2(\{R\leq |x|\leq R+1\}\times [0,1])}\leq
Ne^{-a\frac{R^2}{36d}}.
$$
\end{theorem}
\begin{proof}
We choose $R_0\geq 1$ sufficiently large such that the corresponding
norms in $\{x:|x|\geq R_0\}\times [0,1]$ of $V_1,V_2,p_1,p_2$ are
smaller than $1/(4C)$  where $C$ is the constant in Theorem
\ref{thm1}. Let $\mu:\bR^d\to [0,1]$ be a smooth function such that
$\mu\equiv 0$ for $|x|\leq 1$ and $\mu\equiv 1$ for $|x|\geq 2$. For
any $R>0$, we set
$$\mu_R(x)=\mu(x/R),\quad u_R=u\mu_R,$$
$$V_{R,j}=V_j\chi_{|x|\geq R},
\quad p_{R,j}=p_j\chi_{|x|\geq R},\,\,j=1,2.$$ A direct computation
shows that $H_{V_R,p_R}u_R=e_R$, where
$$
e_R=-\nabla\mu_R\cdot p_1u-\nabla\mu_R\cdot p_2\bar u+\Delta
\mu_Ru+2\nabla \mu_R\cdot\nabla u.
$$
We set $\beta=aR/3,\lambda=C_0\beta^m$, and use Corollary \ref{cor1}
with $u_R,V_R,p_R$ in place of $v,V$ and $p$ and obtain
\begin{equation}
                             \label{eq7.31}
\|e^{\beta\phi_{\lambda}(|x|/\sqrt{d})}u_R\|_{L^2_xL^2_t}+
\|e^{\beta\phi_{\lambda}(|x|/\sqrt{d})}|\nabla
u_R|\|_{B^{\infty,2}_xL^2_t}\leq N(I_1+I_2),
\end{equation}
where
$$
I_1=\|e^{\beta\phi_{\lambda}(|x|)}e_R\|_X,\quad I_2=
\lambda^3\sum_{j=0}^1\|J^1(e^{\beta\phi_{\lambda}(|x|)}u_R(x,j))\|_{L_x^2}.
$$

{\it Estimate of $I_2$:} It is clear that
$$
I_2\leq \lambda^3(\beta+1)\big(\sum_{j=0}^1\|e^{\beta
\phi_{\lambda}(|x|)}u_R(x,j)\|_{L_x^2}
+\|e^{\beta\phi_{\lambda}(|x|)}|\nabla u_R (x,j)|\|_{L_x^2}\big)
$$
On the set $\{x:|x|\geq R\}$, we have
$$\beta\phi_{\lambda}(|x|)\leq
a|x|^2/3\leq a|x|^2-2aR^2/3.$$ Therefore, for $R$ sufficiently large
$$
I_2\leq N(d)
a^{m+1}R^{m+1}e^{-2aR^2/3}\sum_{j=0}^1\|u(x,j)\|_{H^1(e^{a|x|^2})}\leq
1.
$$ This gives the estimate of $I_2$.

{\it Estimate of $I_1$:} Since $e_R$ is supported in  the annulus
$\{R\leq |x|\leq 2R\}$, H\"older's inequality yields
$$
I_1\leq e^{2aR^2/3}\|-\nabla\mu_R\cdot p_1u-\nabla\mu_R\cdot p_2\bar
u+\Delta \mu_Ru+2\nabla \mu_R\cdot\nabla u\|_{B^{1,2}_xL^2_t}
$$
$$
\leq e^{2aR^2/3}\big[\|u\|_{L_x^2L_t^2}\big(\||\nabla
\mu_R|(|p_1|+|p_2|)\|_{B_x^{2,\infty}L_t^\infty}
+\|\Delta\mu_R\|_{B_x^{2,\infty}}\big)
$$
$$
+2\|\nabla u\|_{L^{2}_xL^2_t}\|\nabla \mu_R\|_{B^{2,\infty}_x}\big]
$$
$$
\leq N(d)e^{2aR^2/3}R^{d/2-1}(1+\||p_1|\|_{L^\infty_{t,x}}+
\||p_2|\|_{L^\infty_{t,x}})\|u\|_{H^1_xL_t^2}\leq e^{aR^2},
$$
for $R$ large enough. This completes the estimate of $I_1$, and
using \eqref{eq7.31} we get for $R$ sufficiently large
\begin{equation}
                                \label{eq8.23}
\|e^{\beta\phi_{\lambda}(|x|/\sqrt{d})}u_R\|_{L^2_xL^2_t}+
\|e^{\beta\phi_{\lambda}(|x|/\sqrt{d})}|\nabla
u_R|\|_{B^{\infty,2}_xL^2_t}\leq N(d)e^{aR^2}.
\end{equation}

Notice that in the annulus $\{x;\,6\sqrt{d}R\leq |x|\leq
6\sqrt{d}R+1\}$, we have $u_R=u$ and
$\beta\phi_{\lambda}(|x|/\sqrt{d})\geq 2aR^2$. From \eqref{eq8.23},
we deduce
$$
e^{2aR^2}\|u_R\|_{H^1(\{6\sqrt{d}R\leq |x|\leq
6\sqrt{d}R+1\})L^2_t}\leq N(d)e^{aR^2}.
$$
Upon a change of variable, the theorem is proved.
\end{proof}

\mysection{Lower bounds}

In this section, we shall give lower bounds for the $H^1$ norms of
the solutions to the equation \eqref{eq10.21} in the annulus
$\{R\leq |x|\leq R+1\}\times [0,1]$. We follow the line in Section 3
of \cite{EKPV2}. The following lemma is similar to Lemma 3.1  of
\cite{EKPV2}.
\begin{lemma}
                    \label{lemma3.1}
Assume that $R>0$ and $\varphi:[0,1]\to \bR$ is a smooth function.
Then there exist
$c=c(d,\|\varphi'\|_{L^\infty},\|\varphi''\|_{L^\infty})>0, N=N(d)$
such that the estimate
$$
\alpha^{3/2}R^{-2}\||xR^{-1}+\varphi(t)e_1|e^{\alpha|xR^{-1}+\varphi(t)e_1|^2}v\|_{L_x^2L_t^2}
$$
\begin{equation}
                    \label{eq1.43}
+\alpha^{1/2}R^{-1}\|e^{\alpha|xR^{-1}+\varphi(t)e_1|^2}\nabla
v\|_{L_x^2L_t^2} \leq
N\|e^{\alpha|xR^{-1}+\varphi(t)e_1|^2}Hv\|_{L_x^2L_t^2}
\end{equation}
holds for any $\alpha\geq cR^2$ and $v\in C_0^\infty(\bR^{n+1})$
supported in the set
$$
\{(x,t):\vert xR^{-1}+\varphi(t)e_1\vert \geq 1\}.
$$
\end{lemma}
The proof of this lemma essentially follows that of Lemma 3.1 of
\cite{EKPV2}. However, the second term on the left-hand side of
\eqref{eq1.43} doesn't appear in that Lemma. For the sake of
completeness, we provide the proof.
\begin{proof}
Set $f=e^{\alpha|xR^{-1}+\varphi(t)e_1|^2}v$. Direct computation
yields
$$
e^{\alpha|xR^{-1}+\varphi(t)e_1|^2}Hv=S_\alpha f-4\alpha A_\alpha f,
$$
where
$$
S_\alpha=H+4\alpha^2R^{-2}|xR^{-1}+\varphi(t)e_1|^2,
$$
$$
A_\alpha=R^{-1}(xR^{-1}+\varphi(t)e_1)\cdot \nabla+dR^{-2}/2
+i\varphi'(x_1R^{-1}+\varphi)/2,
$$
satisfying
$$
A_{\alpha}^*=-A_\alpha,\quad S_{\alpha}^*=S_\alpha.
$$
Therefore
$$
\|e^{\alpha|xR^{-1}+\varphi(t)e_1|^2}Hv\|_{L_x^2L_t^2}^2
=\left\langle \,S_{\alpha}f-4\alpha A_{\alpha}f,S_{\alpha}f-4\alpha
A_{\alpha}f\right\rangle
$$
\begin{equation}
                    \label{eq2.20}
\geq -4\alpha\left\langle (S_\alpha A_\alpha-A_\alpha
S_\alpha)f,f\right\rangle =-4\alpha\left\langle
[S_\alpha,A_\alpha]f,f\right\rangle.
\end{equation}
Since
$$
[S_\alpha,A_\alpha]=2R^{-2}\Delta-8\alpha^2 R^{-4}|xR^{-1}+\varphi
e_1|^2
$$
$$
-[(x_1R^{-1}+\varphi)\varphi''+(\varphi')^2]/2+2i\varphi'R^{-1}\partial_{x_1},
$$
the left-hand side of \eqref{eq2.20} is greater or equal to
$$
32\alpha^3R^{-4}\int |xR^{-1}+\varphi(t)e_1|^2|f|^2\,dxdt +8\alpha
R^{-2}\int |\nabla f|^2\,dxdt
$$
$$
+2\alpha \int[(x_1R^{-1}+\varphi)\varphi''+(\varphi')^2]|f|^2\,dxdt
-8\alpha iR^{-1}\int \varphi'\partial_{x_1}f\bar f\,dxdt.
$$
$$
\geq 32\alpha^3R^{-4}\int |xR^{-1}+\varphi(t)e_1|^2|f|^2\,dxdt
+8\alpha R^{-2}\int |\nabla f|^2\,dxdt
$$
$$
-N\alpha\int (1+|x_1R^{-1}+\varphi|) |f|^2\,dxdt-N\alpha R^{-1}\int
|\nabla f||f|\,dxdt,
$$
$$
\geq 16\alpha^3R^{-4}\int |xR^{-1}+\varphi(t)e_1|^2|f|^2\,dxdt +
4\alpha R^{-2}\int |\nabla f|^2\,dxdt,
$$
for $\alpha\geq cR^2$, with
$N=N(\|\varphi'\|_{L^\infty},\|\varphi''\|_{L^\infty})$. Note that
the last inequality is a consequence of the support property of $v$
and the Cauchy-Schwarz inequality. Using the inequality
$|x-y|^2+y^2\geq \epsilon|x|^2+\epsilon|y|^2,$ for some
$\epsilon=\epsilon(d)>0$, the estimate \eqref{eq1.43} follows.
\end{proof}

We are now ready to obtain the following lower bound.

\begin{theorem}
                        \label{thm3.2}
Let $u\in C([0,1]:H^1)$ be a solution of \eqref{eq10.21}. If for
some $\delta_1,\delta_2,L_1,L_2>0$ and $R_1\geq 1$
\begin{equation}
                    \label{eq5.5.10}
\int_{2\delta_2}^{1-2\delta_2}\int_{|x|\leq R_1}|u|^2\,dxdt\geq
\delta_1,
\end{equation}
$$
\|V_1\|_{L^{\infty}_{t,x}}+\|V_2\|_{L^{\infty}_{t,x}}\leq L_1,\quad
\||p_1|\|_{L^{\infty}_{t,x}}+\||p_2|\|_{L^{\infty}_{t,x}}\leq L_2
$$
then there exist $R_0\geq R_1+1$ depending on
$d,\delta_1,\delta_2$,$L_1$ and $A:=\|u\|_{L^2([0,1]:H^1)}$, and a
constant $c_1\geq 1$ depending only on $d$ and $L_2$ such that for
any $R\geq R_0$
\begin{equation}
\delta(R):=\|u\|_{L^2([0,1]:H^1(\{R-1\leq|x|\leq R\}))}\geq
e^{-c_1R^2}.
\end{equation}
\end{theorem}
\begin{proof}
Following \cite{EKPV2}, for any $R>1$ we define smooth functions
$\theta_R,\eta\in C_0^\infty(\bR^n)$ and $\varphi\in
C_0^{\infty}([0,1])$ satisfying
$$\left\{
  \begin{array}{ll}
    \theta_R=1, & \text{if}\,\,|x|\leq R-1\\
    \theta_R=0, & \text{if}\,\,|x|\geq R
  \end{array}
\right.,\quad \left\{
  \begin{array}{ll}
    \eta=0, & \text{if}\,\,|x|\leq 1\\
    \eta=1, & \text{if}\,\,|x|\geq 2
  \end{array}
\right.,
$$
$$
0\leq \varphi\leq 3,\quad |\phi'|\leq 4\delta_2^{-1},\quad \left\{
  \begin{array}{ll}
    \varphi=0, & \text{in}\,\,[0,\delta_2]\cup[1-\delta_2,1]\\
    \varphi=3, & \text{if}\,\,[2\delta_2,1-2\delta_2]
  \end{array}
\right..
$$
We apply Lemma \ref{lemma3.1} to the function
$$
v(x,t)=\theta_R(x)\eta(xR^{-1}+\varphi(t)e_1)u(x,t),\quad (x,t)\in
\bR^d\times [0,1],
$$
and observe that $v$ is supported in $\{|x|\leq R\}\times
[\delta_2,1-\delta_2]$ and satisfies the hypothesis of Lemma
\ref{lemma3.1}. Moreover, on $B_{R-1}\times [2\delta_2,1-2\delta_2]$
we have $v=u$ and $|xR^{-1}+\varphi e_1 |\geq 2$. Thus,
\begin{equation}
                    \label{eq4.36}
\|e^{\alpha|xR^{-1}+\varphi e_1|^2}v\|_{L^2_xL^2_t} \geq
e^{4\alpha}\|u\|_{L^2_x(\{|x|\leq
R_1\})L^2_t([2\delta_2,1-2\delta_2])}\geq e^{4\alpha}\delta_1.
\end{equation}
We compute
$$
H_{V,p}v=\eta(xR^{-1}+\varphi e_1)(2\nabla \theta_R(x)\cdot \nabla
u+u\Delta \theta_R(x)-(ua_1+\bar ua_2)\cdot \nabla \theta_R(x))
$$
$$
+\theta_R(x)[2R^{-1}\nabla\eta(xR^{-1}+\varphi e_1)\cdot \nabla u
+R^{-2}u\Delta \eta(xR^{-1}+\varphi e_1)
$$
$$
+i\varphi'\partial_{x_1}\eta(xR^{-1}+\varphi e_1)u-(ua_1+\bar
ua_2)\cdot R^{-1}\nabla \eta(xR^{-1}+\varphi e_1)]
$$
\begin{equation}
                        \label{eq4.32}
+R^{-1}\nabla \theta_R\cdot\nabla \eta (xR^{-1}+\varphi e_1)u.
\end{equation}
Notice that the first term on the right-hand side of \eqref{eq4.32}
is supported in $B_R\setminus B_{R-1}\times [0,1]$ where
$|xR^{-1}+\varphi e_1|\leq 4$, while the second and the third terms
are supported in $\{(x,t):1\leq|xR^{-1}+\varphi e_1|\leq 2\}$.
Therefore Lemma \ref{lemma3.1} with $\alpha=cR^2$ yields
$$
N^{-1} (c^{3/2}R\|e^{\alpha|xR^{-1}+\varphi e_1|^2}v\|_{L^2_xL^2_t}+
c^{1/2}\|e^{\alpha|xR^{-1}+\varphi e_1|^2}|\nabla v|\|_{L^2_xL^2_t})
$$

$$
\leq L_1\|e^{\alpha |xR^{-1}+\varphi e_1|^2}v\|_{L^2_xL^2_t}
+L_2\|e^{\alpha|xR^{-1}+\varphi e_1|^2}|\nabla v|\|_{L^2_xL^2_t}
$$
\begin{equation}
                    \label{eq4.53}
+e^{16\alpha}\delta(R)+e^{4\alpha}A\delta_2^{-1}.
\end{equation}
Taking $R\geq R_0(L_1,d)$ and $c\geq c_0(L_2,d)$ sufficiently large,
the first two terms on the right-hand side of \eqref{eq4.53} can be
absorbed in the left-hand side. This combined with \eqref{eq4.36}
yields
$$
N^{-1}c^{1/2}R\delta_1\leq e^{12\alpha}\delta(R)+A\delta_2^{-1}.
$$
To complete the prove of the theorem, it suffices to set $c_1=12c$
and choose  $R\geq R_0(d,L_1,\delta_1,\delta_2,A)$.
\end{proof}

\mysection{Proof of theorem \ref{thm1} and \ref{thm2}.}

With the upper and lower bounds established in Sections 2 and 3, we
are ready to prove Theorem \ref{thm1} and \ref{thm2}.

{\it Proof of Theorem $\ref{thm1}:$} Without loss of generality, we
may assume that $u$ satisfies \eqref{eq5.5.10} for some
$R_1,\delta_1$ and $\delta_2$. Indeed, if \eqref{eq5.5.10} doesn't
hold for any $R_1,\delta_1$ and $\delta_2$, then $u\equiv 0$ on
$\bR^d\times [0,1]$ and there is nothing to prove. By Theorem
\ref{thm2.4}  and \ref{thm3.2}, we have for $R$ sufficiently large
$$
e^{c_1(d,L_{2})(R+1)^2}\leq N(d)e^{-aR^2/36d}.
$$
But this is impossible if $a\geq c_0:=37c_1 d$. Hence the theorem
follows.

{\it Proof of Theorem $\ref{thm2}:$} By using Theorem \ref{thm1} and
considering the difference $u:=u_1-u_2$ of two different solutions
of \eqref{eq10.37}, the proof of the theorem follows that of
Corollary 1.2 of \cite{IonKen2} almost word by word. Therefore, it
will be omitted.

\end{document}